\documentclass[12pt]{article}

      \usepackage{latexsym}
         \usepackage[reqno, namelimits, sumlimits]{amsmath}
         \usepackage{amssymb, amsfonts}
         \usepackage{amsthm}
 \newcommand{\beqn}{\begin{eqnarray}}
 \newcommand{\eeqn}{\end{eqnarray}}
 \newcommand{\be}{\begin{equation}}
 \newcommand{\ee}{\end{equation}}
 \newcommand{\ba}{\begin{array}}
 \newcommand{\ea}{\end{array}}
 \newcommand{\pa}{\partial}
 
 \newcommand{\ci}{\cite}
 \newcommand{\la}{\label}

\newcommand{\Om}{\Omega}
\newcommand{\lb}{\lambda}

\newcommand{\na}{\nabla}
\newcommand{\om}{\omega}

\newcommand{\al}{\alpha}

 \newcommand{\De}{\Delta}

\newcommand{\rb}{\mathbb{R}}

\newcommand{\M}{\mathcal{M}}

\def\R{{\rm I\kern-.1567em R}}
\def\M{{\rm I\kern-.1567em M}}

\def\div {{\rm div\,}}

\def\ess{{\rm ess}}
\def\supp{{\rm supp}}
\def\qfive{\int\limits_{Q(z_0,R)}|v|^5\,dz}
\def\qthreeinf{\ess \hskip-10pt\sup\limits_{t\in ]t_0-R^2,t_0[} \int\limits_{B(x_0,R)}|v(x,t)|^3\,dx}
\def\qthree{R^{-2}\int\limits_{Q(z_0,R)}|v|^3\,dz}
\def\qtwo{R^{-3}\int\limits_{Q(z_0,R)}|v|^2\,dz}
\def\qgrad{R^{-1}\int\limits_{Q(z_0,R)}|\nabla v|^2\,dz}
\def\qinft{\ess\hskip-10pt \sup\limits_{(x,t)\in Q(z_0,R)} \sqrt{t_0-t}\,\,|v(x,t)|}
\def\qinfxt{\ess\hskip-10pt\sup\limits_{(x,t)\in Q(z_0,R)} |x-x_0|\,\,|v(x,t)|}

\newcommand{\enorm}[1]{[\hskip-2.6pt|#1|\hskip-2.6pt ]}
\newcommand{\norm}[3]{\|#1\|^{#2}_{#3}}
\newcommand{\eenorm}[3]{\enorm{#1}^{#2}_{#3}}

 \newtheorem{theorem}{Theorem}[section]
 
 \newtheorem{definition}[theorem]{Definition}
 
 \newtheorem{lemma}[theorem]{Lemma}
 
 \newtheorem{remark}[theorem]{Remark}

 \newtheorem{pro}[theorem]{Proposition}

\title{
On Type I singularities of the local axi-symmetric solutions of the Navier-Stokes equations}

 \author{%G. Koch\footnote{University of Chicago}\and
  G. Seregin\footnote{Oxford University. Supported in part by  the Alexander von Humboldt Foundation}
\and V. \v Sver\'ak\footnote{University of Minnesota. Supported in part by NSF
  Grant DMS-0457061}}

\begin{document}
\maketitle
%\begin{center} {\Large\bf  Local regularity of axially symmetric solutions to the %Navier-Stokes equations}\\
%  \vspace{1cm}
% {\large
% G. Koch, N.  N. Nadirashvili,  G. Seregin, and  V. Sverak }
% \end{center}

  \vspace{1cm}
 \noindent
 {\bf Abstract } Local regularity of axially symmetric solutions to the Navier-Stokes equations is studied. It is shown
 that under certain natural assumptions there are no singularities of Type I.

 \vspace {1cm}

\noindent {\bf 1991 Mathematical subject classification (Amer. Math.
Soc.)}: 35K, 76D.


%\noindent\textbf{Acknowledgement}
%supported by the Alexander von Humboldt Foundation and by  the
%RFFI grant 05-01-00941-a.

\setcounter{equation}{0}
\section{Introduction  }

In this paper we will consider local regularity properties of axi-symmetric solutions
of the 3D Navier-Stokes equations
\be
\label{nse1}
\begin{array}{rcl}
\pa_tv+v\cdot\nabla v +\nabla q -\Delta v & = & 0 \\
\div v & = & 0\,\,.
\end{array}
\ee

Most of the known regularity theory for these equations (and, in fact,
for many other equations) is based on optimal estimates for the linear
part and on treating the nonlinearity as a perturbation which
is (locally) small in a suitable sense. An important role in formulating
suitable smallness conditions is played by certain (local) scale-invariant
quantities. These are the quantities which are invariant under the
scaling symmetry $v(x,t)$, $q(x,t)  \to \lambda v(\lambda x,\lambda ^2 t)$,
$\lambda^2 q(\lambda x, \lambda^2 t)$. The reason why the regularity criteria
should be formulated in terms of the scale-invariant quantities is simple:
The class of regular solutions is invariant under the scaling and
therefore sufficient conditions for membership in this class
should ideally be also invariant under the scaling, or at least they should
scale in the correct way, in that the quantity controlling regularity should
not decrease if we scale the solution with $\lambda >> 1$.

To write down examples of the scale-invariant quantities we recall the following
standard notation. The points of the space-time $\mathbb R^n\times \mathbb R$ will be denoted by
$z=(x,t)$. For $x_0\in \rb^n$ we denote by $B(x_0,R)$ the ball $\{x:\,|x-x_0|<R\}$
and for $z_0=(x_0,t_0)\in\rb ^n\times\rb $ we denote by $Q(z_0,R)$ the parabolic
ball $B(x_0,R)\times ]t_0-R^2,t_0[$. Here are some examples of the scale-invariant
quantities for $n=3$:
\begin{eqnarray}
 \qfive\,\,, \label{qfive}\\
 \qthreeinf \,\,,\label{qthreeinft}\\
\qthree\,\,,\label{qthree}\\
\qtwo\,\,,\label{qtwo}\\
\qgrad\,\,,\label{qgrad}\\
\qinft\,\,,\label{qinft}\\
\qinfxt\,\,.\label{qinfxt}
\end{eqnarray}

 A typical local regularity result says that, under some natural technical assumptions
 \footnote{Such assumption must include some control of pressure, such as $q\in L_{3\over 2}$. Fortunately,
  such control is
  available from energy estimates in most situations.}, a point $z_0$ is a regular point of the solution
  $v$ if a suitable scale invariant quantity $X(z_0,R;v)$ of the type in the examples
  above is sufficiently small for all $R\in]0,R_0[$. In fact, $X$ can be any of the
  quantities above with the exception of (\ref{qtwo}), in which case the validity of the corresponding
  result is open. \footnote{The reason for the difficulties with (\ref{qtwo}) is that the space-time $L_2$ norm of $v$ is not sufficiently
  strong to control the energy flux (unless one can come up with some surprising new
  property of the equations).  Roughly speaking, the energy flux is controlled by the $L_3$
  norm of $v$. Since the energy estimate gives the control of the $L_{10\over3}$ norm of $v$,
  there is some gain in regularity and it is natural to try to bootstrap it and  try control the
  energy flux by some $L_r$ norm with $r<3$. This does work, but the borderline exponent
  $r$ for this argument is $r=5/2$, still quite far from $r=2$ which would be needed
  for a local regularity result with $\qtwo$. See for example \ci{S10}.}

  At the time of this writing, there is no known scale-invariant quantity for which an a-priori
  estimate would be known for general 3D solutions. In fact, all the known estimates can be traced
  back to the energy estimate, which gives bounds in quantities such as
  $\int\limits_{B(x_0,R)}|v(x,t)|^2\,\,dx$ or $\int\limits_{Q(z_0,R)}|\nabla v|^2\,\,dz$, which
  do not have the scaling needed for the existing local regularity theory.
  This is often quoted as the main stumbling block in our understanding
  of the Navier-Stokes regularity. This statement is probably correct, at least as
  a first approximation. However, even if we {\it assume} that scale-invariant estimates
  of natural quantities are available, in many cases we are still unable to prove regularity,
  unless an additional  smallness condition is imposed. For the quantities of the type
  $X(z_0,R;v)$  listed above one can show that $X(z_0,R;v)\le C$ for some $C>0$ (not necessarily small)
  implies regularity for (\ref{qfive}) and (\ref{qthreeinft}), but in the remaining cases
  the known theory requires an additional smallness condition (and, as remarked above, the
  situation with (\ref{qtwo}) is even worse). Moreover, even the proofs of the cases (\ref{qfive})
  and (\ref{qthreeinft}) rest on the fact that the assumptions imply that a certain
  quantity becomes small.\footnote{The finiteness of (\ref{qfive}) implies
  $\lim\limits_{R\to0} \qfive = 0$, which  gives us a small quantity. The finiteness of (\ref{qthreeinft})
  implies (for the solutions of the equation) \hbox{$\lim\limits_{R\to 0} \int\limits_{B(x_0,R)}|v(x,t_0)|^3\,\,dx = 0$}.
  This again gives a small quantity, but in this case it is not easy to exploit it, since
  we essentially have to show that some regularity propagates backwards in time.}

  In this paper we study local regularity results for axi-symmetric solutions of the 3D Navier-Stokes
  under an assumption that a weakened version of quantity (\ref{qinft}) or, respectively, (\ref{qinfxt})
  is finite (but not necessarily small).
  These studies can be thought of as a continuation of the work started in~\cite{CSTY},
  \cite{KNSS}, and \cite{CSTY2}.
The exact assumption which we will use to replace  (\ref{qinfxt}), in the
  axi-symmetric situation, with the $x_3$-axis as the axis of symmetry, is
  \be
  \label{assumptionx}
  \ess\hskip-10pt\sup\limits_{(x,t)\in Q(z_0,R) } \sqrt{x_1^2+x_2^2}\,\, |\bar v(x,t)|\,\,<+\infty
  \ee
  for some $R>0$,
  where $z_0$ lies on the $x_3$-axis and we denote by $\bar v(x,t)$ the projection of the velocity vector $v(x,t)$ into the
  plane passing through $x$ and the axis of symmetry $x_3$.
  Similarly, the exact assumption which will replace (\ref{qinft}) in the axi-symmetric situation,
  with the $x_3$-axis as the axis of symmetry, is
  \be
  \label{assumptiont}
  \ess\hskip-10pt\sup\limits_{(x,t)\in Q(z_0,R)} \sqrt{t_0-t}\,\,|\bar v(x,t)|<+\infty
  \ee
  for some $R>0$,
  where $z_0$ and $\bar v$ are as above.
   Our main results are as follows.

   \begin{theorem}\la{irt1} Assume that  $v\in L_3(Q(z_0,R))$
is an axially symmetric weak solution to the Navier-Stokes equations in $Q(z_0,R)$ such that there exists an associated
pressure field $q\in L_{\frac32}(Q(z_0,R))$. If, in addition,
$v$ satisfies (\ref{assumptiont}), then $z_0$ is a regular point of $v$.
\end{theorem}

\begin{theorem}\la{irt2} Assume that  $v\in L_3(Q(z_0,R))$
is an axially symmetric weak solution to the Navier-Stokes equations in $Q(z_0,R)$ such that there exists an associated
pressure field $q\in L_{\frac32}(Q(z_0,R))$. Suppose that $v$ is essentially bounded in the space-time cylinders of the from $B(x_0,R)\times ]t_0-R^2,t'[$  for each $t'<t_0$, where the bound may depend on $t'$. If, in addition,
$v$ satisfies (\ref{assumptionx}), then $z_0$ is a regular point of $v$.
\end{theorem}

These are local versions of the main results in the paper \cite{KNSS}. Similar (but not identical) results
also appeared in \cite{CSTY} and \cite{CSTY2}.

For completeness, we formulate another theorem, which  is a local version of the corresponding global regularity result in \ci {L} and \ci{Yud}.

\begin{theorem}\la{irt3} Assume that  $v\in L_3(Q(z_0,R))$
is an axially symmetric weak solution to the Navier-Stokes equations in $Q(z_0,R)$ such that there exists an associated
pressure field $q\in L_{\frac32}(Q(z_0,R))$. Suppose that $v$ is essentially bounded in the space-time cylinders of the from $B(x_0,R)\times ]t_0-R^2,t'[$  for each $t'<t_0$, where the bound may depend on $t'$. If, in addition, the field $v$ has no swirl, i.\ e.\
$v=\bar v$, then $z_0$ is a regular point of $v$.
\end{theorem}

  On a conceptual level our method will be close to the one used in \cite{KNSS}, and will rely on the Liouville-type theorems
  established in that paper. However, certain important technical parts will be treated in a different way.

  We first recall some terminology related to the Liouville-type results for the Navier-Stokes
  proved in~\cite{KNSS}. An {\it ancient solution} of the Navier-Stokes equation is a solution
  defined in $\rb^n\times]-\infty,0[$. We are interested in ancient solutions with bounded velocity, see Definition \ref{bd1}. Non-zero solutions of this form
  can be generated by a natural re-scaling and limiting  procedures at a potential singularity, see Section 2.
  The definition of the ancient solutions still allows for the ``parasitic solutions" of the form $u(x,t)=b(t)$ (for any bounded $b\colon ]-\infty,0[\to\rb^n$),
  with the corresponding pressure $p$ given by $p(x,t)=-b'(t)\cdot x$, see Remark \ref{br3'}. To exclude these solutions
  (which - under some natural assumptions - cannot arise from the re-scaling procedures, see Theorem \ref{rt8}) we introduce
  the notion of the ancient mild solutions. These are the ancient solutions which satisfy the natural representation formula
  \be\la{mild}
  u(t)=S(t-t_0)u(t_0)+\int_{t_0}^{t} (S(t-s)\,P \div\,(u(s)\otimes u(s))\,ds\,\,
   \ee
   for some sequence of times $t_0\to-\infty$, where $S$ is the solution operator for the heat equation
   and $P$ is the Helmholtz projection onto the div-free fields. (We remark that the usual integration by parts
   shows that the integral on the left-hand side is well defined for $u\in L_\infty$.) See~\cite{KNSS}
   for details. The strongest conjecture regarding the Liouville-type results one can make about
   the Navier-Stokes equations is the following:

   \bigskip
   \noindent
   Conjecture (L): {\sl The velocity field of any bounded mild ancient solution of the Navier-Stokes equations is constant.}
   \bigskip

   The conjecture was proved for $n=2$ and also for axi-symmetric solutions in 3D, provided the additional
   decay condition
   \be\label{decay}
   \sqrt{x_1^2+x_2^2}\,|v(x,t)|\le C \quad \mbox{in $\rb^3\times]-\infty,0[$}
   \ee
   is satisfied, see \cite{KNSS}.

   What would be the implications of the validity of Conjecture (L) for the regularity theory?
   Roughly speaking, if Conjecture (L) is valid, then all the problems discussed above concerning regularity
   in the presence of a scale-invariant estimates are solved. Indeed, the re-scaling procedure
   preserves any scale-invariant estimate, and typically the estimate will also be preserved
   in the limiting process. Therefore as a result of the re-scaling we get, in the limit,
   a non-zero bounded mild ancient solution for which a scale invariant quantity is finite. Conjecture (L)
   would leave only one candidate for the mild ancient solution - namely a non-zero constant
   velocity field. However, this possibility is typically not compatible with a finite
   scale-invariant bound.

   We can summarize the above as follows:
   \begin{equation*}
   \begin{array}{c}
   \mbox{scale invariant} \\
   \mbox{estimate}
   \end{array}
   +
   \,\,\mbox{Conjecture (L)} \,\,\,\,\Rightarrow\,\,\,\,\mbox{regularity}\,.
   \end{equation*}

   Singularities for which some scale-invariant quantity is bounded are often called
   Type I singularities. (The most common definition of Type I singularities
   uses quantity (\ref{qinft}).)

   While we do not really know what the likelihood of Conjecture (L) being true is
   for the general 3D solutions, we are quite confident that the conjecture is indeed true
   for the axi-symmetric solutions. The axi-symmetric case of Conjecture (L) would
   imply much stronger results than Theorems~\ref{irt1} and~\ref{irt2} above.
   However, we have not been able to fully prove Conjecture (L) in the axi-symmetric case so far.

   Our method of proof of Theorems~\ref{irt1} and~\ref{irt2} is can be described as follows.
   Roughly speaking, we will show that, on the solutions of the equations, the assumed
   scale invariant bounds imply that all the other important scale-invariant quantities
   are bounded, and these bounds, together with the known partial regularity theory
   (\cite{CKN, LS, Li}), lead relatively easily to the bounds required by the Liouville theorems
  in \cite{KNSS}.

   The idea that a bound of one scale-invariant quantity should lead (for the solutions
   of the equations) to bounds on other scale invariant quantities is of course not new.
   However, examples from some other elliptic/parabolic PDEs  show that these issues can be subtle. For example,
   in the theory of harmonic mappings or the harmonic map heat flow we do have a scale-invariant
   a-priori bound, which corresponds to a bound of quantity (\ref{qtwo}). However, it is known
   that singularities can still arise, and therefore the bound corresponding to (\ref{qfive})
   (which is known to imply regularity in that situation) cannot be derived from
   (the analogue of) (\ref{qtwo}).

   We now informally explain the main steps of the proof. One is that the swirl component
   of the velocity field, $v_\varphi=v\cdot e_\varphi$
   %, where $e_\varphi=(-x_2,x_1,0)$,
   satisfies a scalar parabolic equation which enables one to gain some regularity.
   To explain this, we need to introduce the following simple
notation. Let  $e_1$, $e_2$, $e_3$ be an orthogonal basis of the
Cartesian coordinates $x_1$, $x_2$, $x_3$ and $e_\varrho$,
$e_\varphi$, $e_3$ be an orthogonal basis  of the cylindrical
coordinates $\varrho$, $\varphi$, $x_3$ chosen so that
$$e_\varrho=\cos \varphi e_1+\sin \varphi e_2, \quad e_\varphi=-\sin\varphi
e_1+\cos\varphi e_2,\quad e_3=e_3.$$ Then, for any vector-valued
field $v$, we have  representations
$$v=v_ie_i=v_1e_1+v_2e_2+v_3e_3=v_\varrho e_\varrho+v_\varphi e_\varphi+v_3e_3.$$
   Next, letting
   %$r=\sqrt{x_1^2+x_2^2}$, and
   $f=\varrho v_\varphi$, we have
   \be
   \label{f}
  \pa_tf +\bar v\cdot\nabla f=\Delta f - {\frac 2\varrho} {{\partial f}\over{\partial \varrho}}\,\,.
  \ee
We would like to prove a $L_{\infty}$-bound on $f$. Such a bound will give us enough
  information about $v-\bar v$ so that, oversimplifying slightly, we can
  replace $\bar v$ by $v$ in our assumptions. The $L_\infty$ bound for (\ref{f})
  does not follow from general parabolic theory, since the general theory requires
  more regularity than we have. However, it is known that if the drift term in equations
  such as (\ref{f}) is div-free, one can prove the $L_\infty$ estimate for $f$ with
  weaker assumptions on the coefficients. See for example \cite{FrehseRuzicka} for the elliptic case
  and \cite{Zhang} for the parabolic case.

   Another important step in the proof is conceptually the same as deriving an estimate
   for the quantities (\ref{qthree})  and (\ref{qgrad}) from the boundedness of (\ref{qinft}).
   This can be done by bootstrapping the energy inequality. This idea was used for example in
   \cite{S9}.  The  technical details
   are somewhat complicated, but the main idea can be explain at a heuristic level as follows.
   To simplify notation, we will use $Q(R)$ for $Q(0,R)=Q((0,0),R)$, $Q$ for $Q(1)$,
   $B(R)$ for $B(0,R)$, and $B$ for $B(1)$.

   We first note that (\ref{qinft}) implies a bound on the (scale-invariant)
   quantity $R^{1-\frac 2l-\frac 3s} \|v \|_{s,l,Q(R)} $ for $l<2$ and $s\ge 1$. Here, $\|\cdot\|_{s,l,Q(R)} $ is the norm of the mixed Lebesgue space $L_{s,l}(Q(R))=L_l(-R^2,0;L_s(B(R)))$.

   Let $\eenorm u {}{Q(R))}$ denote the parabolic energy norm in $Q(R)$, i.\ e.\
$$
\eenorm u 2 {Q(R)} = \textrm {ess sup}_{t\in]-R^2,0[}\norm {u(\cdot,t)} 2 {L_2(B(R))} + \norm{\nabla u} 2 {L_2(Q(R))}\,.
$$

To avoid technicalities, let us pretend that the pressure satisfies $|q|\sim |v|^2$. In reality it is not quite
true and in the rigorous proof one has to deal with this, but the procedure is well understood. Therefore our
simplifying assumption $|q|\sim |v|^2$ is reasonable for the heuristics.
We will now work with $R=1$, but we can scale the calculations to any $R>0$, if we divide
all the involved quantities by the powers of $R$ which make them scale-invariant.

The local energy inequality implies
\begin{equation}
\label{e1}
\eenorm v 2{Q} \lesssim \norm v 3 {L_3(Q(2))} + \norm v 2 {L_3(Q(2))}\,.
\end{equation}

We can now ``bootstrap" this inequality. There are some technical complications coming
from the fact that we have $Q(2)$ on the right-hand side of \eqref{e1} but only $Q$
on the left-hand side. Such problems come up often in local regularity
theory of elliptic and parabolic equations, and it is quite well-understood how to deal with them
by suitable iteration procedures.  Therefore, cheating slightly, we can pretend that
we actually have $Q$ on the right-hand side of inequality~\eqref{e1}:

\begin{equation}
\label{e2}
\eenorm u 2 {Q}\lesssim \norm v 3 {L_3(Q)} + \norm v 2 {L_3(Q)}\,.
\end{equation}

To bootstrap, we  estimate
\begin{equation}
\label{e3}
\norm v {} {L_3(Q)}\lesssim \eenorm v \alpha Q \norm v \beta {s,l, Q}
\end{equation}
for suitable $\alpha, \beta >0$, $\alpha + \beta =1$, and use this in inequality~\eqref{e2}.
We see that when $\alpha<2/3$, we can estimate $\eenorm v {}{Q}$
%$\eenorm v Q(1) $
in terms of the
norm $\norm v
{}{ s,l,Q}
$. (The process also works for $\alpha=2/3$ provided
$\norm v {} {s,l}$ is sufficiently small.)
\smallskip

It remains to determine the correct exponents $\alpha, \beta$ in~\eqref{e3}. Denoting by $2^*$ the Sobolev exponent
of the space $W^{1,2}$ (i.\ e.\ $1/{2^*}=1/2-1/n$), we have by the H\"older inequality
\begin{equation}
\label{e4}
\norm v {} {3,3}\le \norm v {\alpha_1} {2, \infty }
\norm v {\alpha_2} {{2^*},2}
\norm v {\alpha_3} {s,l}\,\,,
\end{equation}
where $\alpha_1,\alpha_2,\alpha_3$ are non-negative numbers satisfying
\be
\label{e5}
\begin{array}{rcl}
\alpha_1+\alpha_2+\alpha_3 & = & 1 \\
\frac {\alpha_1}2 + {{\alpha_2}\over{2^*}} + {{\alpha_3}\over s} & = & 1\over 3 \\
{{\alpha_2}\over 2} + {{\alpha_3}\over l} & = & 1\over 3
\end{array}
\ee
By Sobolev imbedding, we have from~\eqref{e4}
\be
\label{e6}
\norm v {} {L_3(Q)}\lesssim \eenorm v {\alpha_1+\alpha_2} Q \norm v {\alpha_3} { s, l, Q}\,\,,
\ee
and we see that~\eqref{e3} holds true with $\alpha=\alpha_1+\alpha_2$, $\beta=\alpha_3$.
Therefore the set of the parameters $s,l$ for which the iteration procedure works is given by the
condition that equations~\eqref{e5} for $\alpha_1,\alpha_2,\alpha_3$ have a non-negative solution
with $\alpha_3>1/3$. Solving~\eqref{e5} for $n=3$, we obtain
\newcommand{\xl}[1]{{#1}\over {l}}
\newcommand{\xs}[1]{{#1}\over {s}}
\newcommand{\ddd}{{\xs 3} + {\xl 2} - {3\over 2}}

\begin{eqnarray*}
\alpha_1 & = & {{{\xs 1} + {\xl 1} - {2\over 3}}\over{\ddd}}\,\,,\\
\alpha_2 & = & {{{\xs 2} + {\xl 1} - 1}\over{\ddd}}\,\,,\\
\alpha_3 & = & {{1}\over{6(\ddd)}}\,\,.
\end{eqnarray*}
One can check easily that the conditions $\alpha_1\ge0,\,\alpha_2\ge0, \,\alpha_3>1/3$
are equivalent to
\be
\label{e7}
\begin{array}{ccc}
{\xs 1} + {\xl 1} & \ge & {2\over 3}\,\,,\\
{\xs 2} + {\xl 1} & \ge &  1\,\,,\\
{\xs 3} + {\xl 2} & < & 2\,\,.
\end{array}
\ee
In the plane with coordinates $x={\xs 1}$ and $y={\xl 1}$, the last set of
equations describes a thin triangle contained in the first quadrant.
It is easy to see that one can choose a suitable $l<2$ and $s>1$ for which
these conditions are satisfied.

\smallskip
Above we worked with parabolic balls of radius of order 1.
It is clear that the calculations can be ``scaled" to $Q(R)$ and $Q(2R)$ if
we divide the $L_{s,l}$-norms by suitable powers of $R$ to obtain scale invariant quantities.
It is also clear that of the restrictions in \eqref{e7}, only the last one is crucial,
since by using H\"older inequality one can always move to higher exponents, as long as
the power in the correct scaling factor remains positive.

This finishes our explanation of the heuristics behind the second step of the proof.
We did not mention one more complication. Since our main assumption involves
only $\bar v$ and the information about $v-\bar v$ is obtained from equation (\ref{f}),
we have to use different one set of parameters $l,s$ for the $\bar v$ component of $v$ and
a different set for $v - \bar v$.  However, this is a technical issue which does not
change the heuristics.

Once we know that the scaled energy-type quantities are bounded, it is not difficult to
derive the bounds which we need in the version of the Liouville conjecture for
axi-symmetric solutions which was proved in~\cite{KNSS}.

%There is one last technicality we need to deal with, namely that the assumptions in
%Theorems~\ref{irt1}, \ref{irt2}, and \ref{irt3} guarantee that our solutions have %enough regularity
%for the above reasoning to be valid.
%(For example, we need to know that they in fact are the {\it suitable weak %solutions}
%in the sense of \cite{CKN}, \cite{Li}, or \cite{LS}. This is a standard application %of
%the usual bootstrapping argument, modulo the fact that some caution is needed due to %the
%potential presence of the local analogy of the parasitic solutions mentioned above,
%namely the solutions of the form $v(x,t)=\nabla h(x,t)$, $p(x,t)=-\pa_th(x,t)$, %where
%the function $h$ is harmonic in $x$, but can have arbitrary dependence on $p$.
%(These solutions also show that in local regularity theory one cannot avoid %assumptions
%about the pressure.)

 Another aim of the paper is to give an alternative approach to certain technical issues  arising in the study  of bounded ancient solutions. In the approach here, we do not use the exact representation formulae which were used in \ci{KNSS}. It turned out to be quite convenient to describe differentiability properties of bounded ancient solutions in terms of certain ``uniform'' Lebesgue and Sobolev spaces, compare with \ci{LR1}. We hope that both approaches are of interest, and complement each other.

\setcounter{equation}{0}
\section{Preliminaries} In this Section,
we recall  known definitions of (weak) solutions to the Navier-Stokes.
\begin{definition}\la{rd1} A weak solution to the Navier-Stokes equations in a domain $\mathcal O\subset \mathbb R^n\times ]t_1,t_2[$ is a divergence free vector-valued field $v\in L_{2,loc}(\mathcal O)$ satisfying
$$\int\limits_{\mathcal O}(v\cdot \pa_tw+v\otimes v:\na\,w+v\cdot \De w)dx\,dt=0$$ for any solenoidal vector-valued field $w\in C^\infty_0(\mathcal O)$.\end{definition}

An important family of weak solutions is given by $v(x,t)=\na h(x,t)$ where $h:\mathcal O\to\mathbb R$ satisfies $\De h=0$. (Dependence on $t$ can be arbitrary). This example shows that further assumptions are needed to obtain some regularity of solutions in the time direction.

Very often, we shall study local regularity of solutions to the Navier-Stokes equations in the unit parabolic ball $Q=B\times ]-1,0[$, where $B=B(1)=B(0,1)$. It is not a loss of generality because of the Navier-Stokes scaling.

In local analysis, the most reasonable object to study is
so-called suitable weak solutions, introduced by
Caffarelli-Kohn-Nirenberg in \ci{CKN}. We
are going to use a slightly simpler  definition of F.-H. Lin in
\ci{Li}
\begin{definition}\la{asd1} The pair $v$ and $q$ is called a
suitable weak solutions to the Navier-Stokes equations in $Q$ if
the following conditions are satisfied:
$$v\in L_{2,\infty}(Q)\cap W^{1,0}_2(Q),\qquad q \in L_\frac
32(Q);$$
$$v\,\mbox{and}\,q \,\mbox{satisfy the Navier-Stokes equtions in
the sense of distributions};$$ for a.a. $t\in ]-1,0[$, the local
energy inequality
$$\int\limits_\mathcal C\varphi(x,t)|v(x,t)|^2dx+2\int\limits^t_{-1}
\int\limits_\mathcal C\varphi |\na v|^2dxdt'\leq
\int\limits^t_{-1} \int\limits_\mathcal C\Big\{|v|^2(\De \varphi
+\pa_t\varphi)+$$
$$+v\cdot\na \varphi(|v|^2+2q)\Big\}dxdt'$$
holds for all non-negative cut-off functions $\varphi\in
C^\infty_0(\mathbb R^3\times \mathbb R)$ vanishing in a
neighborhood of the parabolic boundary of $Q$.
\end{definition}
Here, the following functional spaces have been used:
$$L_{s,l}(Q)=L_l(-1,0;L_s(B)),\qquad L_s=L_{s,s},$$
$$W^{1,0}_{s,l}(Q)=\{v,\,\na v \in L_{s,l}(Q)\}, \qquad W^{1,0}_{s}=W^{1,0}_{s,s},$$
$$W^{2,1}_{s,l}(Q)=\{v,\,\na v ,\,\pa_tv,\,\na^2v\in L_{s,l}(Q)\},\qquad W^{2,1}_{s}=W^{2,1}_{s,s}.$$ The norm of the space $L_{s,l}(Q)$ is denoted by
$\|\cdot\|_{s,l,Q}$.

For further discussions of Definition \ref{asd1}, we refer the reader to
papers \ci{LS} and \ci{S8}.

In what follows,  we shall assume  $v$ and $q$ satisfy the following standing conditions:
  $$\mbox{pair}\,\,v\in L_3(Q)\,\,\mbox{and}\,\,q\in L_\frac 32(Q)\,\,\mbox{satisfies the Navier-Stokes equations}$$
    \be\la{b8}\mbox{in the sense of distributions};\ee
\be\la{b9}v\in L_\infty(B\times ]-1,-a^2[)\qquad \forall a\in ]0,1[;\ee
\be\la{b10} \mbox{there is a number}\,\, 0<r_1<1\,\,\mbox{such that}\,\,v\in L_\infty(Q_1),\ee where $Q_1=B_1\times ]-1,0[$, $B_1=\{r_1<|x|<1\}$.
To explain why there is no loss of generality, we first notice that the pair $v$ and $q$, satisfying conditions (\ref{b8})--(\ref{b10}), is in fact a suitable weak solution to the Navier-Stokes equations in $Q$. It is certainly true in $B\times ]-1,-a^2[$ but condition (\ref{b8}) allows us to extend this property to the whole cylinder $Q$.

Consider now an arbitrary suitable weak solution $v$ and $q$ in $Q$. Let $S\subset B\times ]-1,0]$ be a set of singular points of $v$. It  is closed in $Q$. As it was shown in \ci{CKN}, ${\mathcal P}^1(S)=0$, where ${\mathcal P}^1$ is the one-dimensional parabolic Hausdorff measure. Let us assume that $S\neq\emptyset$. We can choose number $R_1$ and $R_2$ satisfying  $0<R_1<R_2<1$ such that
 $S\cap\overline{Q(R_1)\setminus Q(R_2)}=\emptyset$ and $S\cap B(R_2)\times ]-R^2_2,0]\neq\emptyset$. We put
 $$t_0=\inf\{t\,\,:\,\,(x,t)\in S\cap B(R_2)\times ]-R^2_2,0]\}.$$
 Clearly, $(x_0,t_0)\in S$ for some $x_0\in B(R_2)$. In a sense, $t_0$ is the first singular time of our suitable weak solution $v$ and $q$ in $Q(R_1)$. Next, the one-dimensional Hausdorff measure of the set $$S_{t_0}=
 \{x_*\in B(R_2)\,\,:\,\, (x_*,t_0) \,\,\mbox{is a singular point}\, \}$$ is zero as well. Therefore, given $x_0\in S_{t_0}$, we can find sufficiently small $0<r < \sqrt{R_2^2+t_0}$ such that $B(x_0,r)\Subset B(R_2)$ and
 $\pa B(x_0,r)\cap S_{t_0}=\emptyset$. Since the velocity field $v$ is H\"older continuous at regular points, we can ensure that all conditions of type  (\ref{b8})--(\ref{b10}) hold in the parabolic ball $Q(z_0,r)$ with $z_0=(x_0,t_0)$. We may shift and re-scale our solution if $x_0=0$ and $r\neq 1$.
  %So, in what follows, without loss of generality,  we may assume  $v$ and $q$ %satisfy the following standing conditions:
  %$$\mbox{pair}\,\,v\in L_3(Q)\,\,\mbox{and}\,\,q\in L_\frac  %32(Q)\,\,\mbox{satisfies the Navier-Stokes equations}$$
  %\be\la{b8}\mbox{in the sense of distributions};\ee
%\be\la{b9}v\in L_\infty(B\times ]-1,-a^2[)\qquad \forall a\in ]0,1[;\ee
%\be\la{b10} \mbox{there is a number}\,\, 0<r_1<1\,\,\mbox{such that}\,\,v\in %L_\infty(Q_1),\ee where $Q_1=B_1\times ]-1,0[$, $B_1=\{r_1<|x|<1\}$.

 In our investigations of regularity of suitable weak solutions, the particular case of weak solutions, see Definition \ref{rd1}, plays a crucial role. Here is the corresponding definition.
 \begin{definition}\la{bd1} (\ci{KNSS}) A bounded divergence free field $u\in L_\infty(Q_-;\mathbb R^n)$ is called a weak bounded ancient solution (or simply bounded ancient solution) to the Navier-Stokes equations if
$$\int\limits_{Q_-}(u\cdot \pa_t w+u\otimes u:\na w+u\cdot \De w)dz=0$$
for any $w\in {\stackrel{\circ}{C}} {^\infty_{0}(Q_-)}$. \end{definition}
Here, we have used the following notation:
$$Q_-\equiv\mathbb R^n\times ]-\infty,0[
,\qquad {\stackrel{\circ}{C}} {^\infty_{0}(Q_-)}=\{v\in C^\infty_0(Q_-),\quad \div v=0\}.$$
%Q^{t_0}\equiv\mathbb R^n\times ]t_0-1,t_0[.$$
%and
%$${\stackrel{\circ}{C}} {^\infty_{0}(\mathbb R^n)}=\{v\in C^\infty_0(\mathbb %R^n),\quad \div v=0\},$$

The notion of bounded ancient solutions is not quite satisfactory for our purposes since it allows ``parasitic solutions" $u(x,t)=b(t)$, where $b$ is an $L_\infty$-function. (These correspond to harmonic function $h(x,t)=b(t)\cdot x.)$

The important subclass of bounded ancient solution was introduced in \ci{KNSS}. It consists of the so-called mild bounded ancient solutions, i.e., bounded ancient solution satisfying the representation formula (\ref{mild}). In \ci{KNSS}, we also showed that one has a natural decomposition $u(x,t)=w(x,t)+b(t)$, where $w$ is given by the right hand side of (\ref{mild}) on a suitable interval $]t_1,t_*[$. (In particular, $w$ is H\"older continuous.)

In this paper, we give another proof of the decomposition of arbitrary bounded ancient solutions into regular and singular parts, see Section 5. It is based on  recovery of a pressure field associated with a given bounded ancient solution. To formulate our theorem about the pressure, we need 
%additional notation
to introduce certain functional spaces

By $L_m(\Om)$ and $W^1_m(\Om)$, we denote the usual Lebesgue and Sobolev spaces of functions defined on $\Om\in \mathbb R^n$. We also need parabolic versions of Sobolev's spaces:
$$W^{1,0}_m(Q_T)=\{|v|+|\na v|\in L_m(Q_T)\},$$$$ W^{2,1}_m(Q_T)=\{|v|+|\na v|
+|\pa_t v|+|\na^2v|\in L_m(Q_T)\}$$
where $Q_T=\Om\times ]0,T[$.
%, $\na $ is the spatial gradient, and $\pa_t$ is the partial derivative in time $t$
The norm of the space $L_m(\Om)$ is denoted by
$\|\cdot\|_{m,\Om}$.

We also going to use the following "uniform" spaces (compare with \ci{LR1}):
%. First, we introduce the following ``uniform" spaces (compare with \ci{LR1}):
$${\cal L}_m(Q_-)=\{\|f\|_{{\cal L}_m(Q_-)}=\sup\limits_{z_0\in Q_-}\|f\|_{m,Q(z_0,1)}<+\infty\},$$
$${\cal W}^{1,0}_m(Q_-)=\{\|f\|_{{\cal W}^{1,0}_m(Q_-)}=\sup\limits_{z_0\in Q_-}\|f\|_{W^{1,0}_m(Q(z_0,1))}<+\infty\},$$
$${\cal W}^{2,1}_m(Q_-)=\{\|f\|_{{\cal W}^{2,1}_m(Q_-)}=\sup\limits_{z_0\in Q_-}\|f\|_{W^{2,1}_m(Q(z_0,1))}<+\infty\}.$$

 To define the regular part of the pressure, we recall the known fact, see e.g. \ci{St}.
 Given $F=L_\infty(\mathbb R^n;\mathbb M^{n\times n})$, there exists a unique function $p_F\in BMO(\mathbb R^n)$ such that $[p_F]_{B(1)}=0$ ($[f]_\Om$
is the mean value of a function $f$ over a spatial domain $\Om\in\mathbb R^n$) and $$\De p_F=-\div \div F=F_{ij,ij} \qquad \mbox{in}\quad \mathbb R^3$$
in the sense of distributions. Moreover, function $p_F$ meets the estimate
$$\|p_F\|_{BMO(\mathbb R^n)}\leq c(n)\|F\|_{\infty,\mathbb R^n}.$$
So, given a bounded ancient solution $u$, we define a regular part of the pressure
 $p_{u\otimes u}$.
\begin{theorem}\la{bt2} Let $u$ be an arbitrary bounded ancient solution. For any number $ m>1$,
$$|\na u|+|\na^2u|+|\na p_{u\otimes u}|\in {\cal L}_m(Q_-).$$
In addition, for each $t_0\leq 0$, there exists a function $b_{t_0}\in L_\infty(t_0-1,t_0)$ with the following property
$$\sup\limits_{t_0\leq 0}\|b_{t_0}\|_{L_\infty(t_0-1,t_0)}\leq c(n)<+\infty.$$

If we let $u^{t_0}(x,t)=u(x,t)+b_{t_0}(t),\,(x,t)\in Q^{t_0}=\mathbb R^n\times ]t_0-1,t_0[$, then, for any number $m>1$,
%the  estimate
$$\sup\limits_{z_0=(x_0,t_0),\,x_0\in \mathbb R^n,\,t_0\leq 0}\|u^{t_0}\|_{W^{2,1}_m(Q(z_0,1))}\leq c(m,n)<+\infty.$$
 %is valid.
  Moreover, for each $t_0\leq 0$, functions $u$ and $u^{t_0}$ obey the system of equations
$$\pa_t u^{t_0}+\div u\otimes u -\De u=-\na p_{u\otimes u},\qquad \div u=0$$ a.e in $ Q^{t_0}$.\end{theorem}
\begin{remark}\la{br3} The first equation of the latter system can be reduced to the form
%rewritten in the following way
$$\pa_t u+\div u\otimes u -\De u=-\na p_{u\otimes u}-b'_{t_0}\quad \mbox{in}\quad Q^{t_0},$$
which is understood in the sense of distributions, $b'_{t_0}(t)=d b_{t_0}(t)/dt$. So, the real pressure field in $Q^{t_0}$ is the  distribution $p_{u\otimes u} +b'_{t_0}\cdot x$.\end{remark}
\begin{remark}\la{br3'}
We can find a measurable vector-valued function $b$ defined on $]-\infty,0[$ and
having the following property.  For any $t_0\leq 0$, there exists a constant vector
$c_{t_0}$ such that
$$\sup\limits_{t_0\leq 0}\|b-c_{t_0}\|_{L_\infty(t_0-1,t_0)}<+\infty.$$
Moreover, the Navier-Stokes system takes the form
$$\pa_t u+\div u\otimes u -\De u=-\na (p_{u\otimes u}+b'\cdot x),\qquad \div u=0$$
in $Q_-$ in the sense of distributions.
\end{remark}
\begin{remark}\la{br4} %In most of  our applications, we shall have some additional
%global information about the pressure field, which will make it possible to conclude
%that $b'=0$.
%all the result of Theorem \ref{bt2} are valid for $b_{t_0}=0$ for any $t_0\leq 0$.
%For example, it is true if the pressure field belongs to $ %L_\infty(-\infty,0;BMO(\mathbb R^n))$.
Bounded ancient solutions with $b'=0$ were introduced in \ci{KNSS} and called mild ancient solutions.  They were systematically studied in \ci{KNSS} and in particular it was shown there that mild ancient solutions are infinitely smooth.\end{remark}

Our interest in bounded ancient solutions comes from their appearance as natural limits of suitable re-scaling procedures at potential singularities. In the context of solutions to the Navier-Stokes equations in $\mathbb R^n\times ]t_1,t_2[$, this was studied in \ci{KNSS}.

We shall now show  that re-scaling procedure also works for potential singularities of local suitable weak solutions. The important point will be that, even in the local (in space) situation, the solutions arising from the re-scaling procedure at a potential singularity are still mild bounded ancient solutions. To be more precise, let us consider local solutions of the Navier-Stokes equations satisfying assumptions
 (\ref{b8})--(\ref{b10}) and introduce functions
$$G(t)=\max\limits_{x\in \overline{B}(r_1)} |v(x,t)|,\qquad M(t)=\sup\limits_{-1\leq \tau\leq t}G(t).$$ Assume that there are singular points of $v$ which are located somewhere on the set $\{ (x,0)\,\,:\,\, |x|\leq r_1\}$. By known regularity criteria (see e.g. \ci{S10}), we have
$$M(t)>\frac \varepsilon{\sqrt {-t}}$$
for some $\varepsilon>0$ and thus
$$M(t)\to +\infty$$
if $t\to 0-$.  We can construct a sequence $t_k$ such that $t_k\in ]-1,0[$, $t_k<t_{k+1}$, $t_k\to 0$, and
$$M(t_k)=G(t_k)=|v(x_k,t_k)|\to +\infty
$$
for some $x_k\in \overline{B}(r_1)$.

Next, we scale   $v$ and $q$ the following way
$$u^{k}(y,s)=\lambda_kv(x,t),\qquad p^k(y,s)=\lambda^2_kq(x,t),$$
where $x=x^k+\lambda_ky$, $t=t_k+\lambda^2_ks$, and $\lambda_k=1/M_k$.
The ball $B(r)$ is mapped  by the change of variables onto $B(-x^k/\lambda_k, r/\lambda_k)$ and if $r\in ]r_2,1]$ then, given any  $R>0$,
$$B(R)\subset B\Big(-\frac {x^k}{\lambda_k}, \frac r{\lambda_k}\Big)$$
for sufficiently large $k$.
For scaled functions $u^k$ and $p^k$, we know
$$ u^k\,\,\mbox{and}\,\,p^k\,\,\mbox{satisfy the Navier-Stokes equations and} $$
\be\la{b15}|u^k|\leq 1\qquad \mbox{in}\,\, B\Big(-\frac {x^k}{\lambda_k}, \frac 1{\lambda_k}\Big)\times \Big]-\frac {1-t_k}{\lambda_k^2},0\Big[;\ee
\be\la{b16} |u^k(0,0)|=1.\ee
\begin{theorem}\la{rt8} For each $ a>0$, the sequence $\{u^{k}\} $ is uniformly H\"older continuous on the closure of $Q(a)$ for sufficiently large $k$ and a subsequence $\{u^{k_j}\} $ of $\{u^{k}\} $ converges uniformly on compact subsets of $\mathbb R^n\times ]-\infty,0]$ to a mild bounded ancient solution $u$ with $|u(0,0)|=1$.\end{theorem}
\begin{remark} The main point of the theorem is that the limit $u$ is a mild solution, i.e., the parasitic solutions cannot appear in this re-scaling procedure. \la{rr9}\end{remark}
\textsc{Proof of Theorem \ref{rt8}}
Our solution $v$ and $q$ has good properties inside $Q_1$. Let us enumerate them.  Let $Q_2=B_2\times ]-\tau^2_2,0[$, where $0<\tau_2<1$, $B_2=\{r_1<r_2<|x|<a_2<1\}$. Then, for any natural $k$,
$$z=(x,t)\mapsto \na^kv(z)\,\,\mbox{is H\"older continuous in}\,\,\overline{Q}_2;$$
$$q\in L_\frac 32(-\tau^2_2,0;C^k(\overline{Q}_2)).$$
The corresponding norms are estimated by constants depending on $\|v\|_{3,Q}$, $\|q\|_{\frac 32,Q}$, $\|v\|_{\infty,Q_1}$, and numbers $k$, $r_1$, $r_2$, $a_2$, $\tau_2$. In particular, we have
\be\la{b11} \max\limits_{x\in \overline{B}_2}\int\limits_{-\tau^2_2}^0|\na q(x,t)|^\frac 32 dt \leq c_1<\infty.\ee
Proof of this statements can be done by induction and founded in \ci{ESS4}, \ci{LS},
and \ci{NRS}.

Now, let us decompose the pressure $q=q_1+q_2$. For $q_1$, we have
$$\De q_1(x,t)=-\div\div \Big[\chi_B(x) v(x,t)\otimes v(x,t)\Big],\quad x\in\mathbb R^3,\quad -1<\tau<0,$$
where  $\chi_B(x)=1$ if $x\in B$ and $\chi_B(x)=0$ if $x\notin B$. Obviously, the estimate
$$\int\limits^0_{-1}\int\limits_{\mathbb R^3}|q_1(x,t)|^\frac 32dxdt\leq c\int\limits_Q|v|^3dz$$ holds which is a starting point for local regularity of $q_1$. Using essentially the same bootstrap arguments, we can show
\be\la{b12}\max\limits_{x\in \overline{B}_3}\int\limits_{-\tau^2_2}^0|\na q_1(x,t)|^\frac 32 dt \leq c_2<\infty,\ee where $B=\{r_2<r_3<|x|<a_3<a_2\}$.
From (\ref{b11}) and (\ref{b12}), it follows that
\be\la{b13}\max\limits_{x\in \overline{B}_3}\int\limits_{-\tau^2_2}^0|\na q_2(x,t)|^\frac 32 dt \leq c_3<\infty.\ee But clearly $q_2$ is a harmonic function in $B$, thus, by the maximum principle,
we have
\be\la{b14}\max\limits_{x\in \overline{B}(r_4)}\int\limits_{-\tau^2_2}^0|\na q_2(x,t)|^\frac 32 dt \leq c_3<\infty,\ee where $r_4=(r_3+a_3)/2$.

Let us re-scale each part of the pressure separately, i.e.,
$$p_{i}^k(y,s)=\lambda_k^2q_i(x,t),\quad i=1,2,$$
so that $p^k=p^k_1+p^k_2$. As it follows from (\ref{b14}), for $p^k_2$, we have
\be\la{b17}\sup\limits_{y\in B(-x^k/\lambda_k, r_4/\lambda_k)}\int\limits^0_{-(\tau^2_2-t_k)/\lambda^2_k}|\na_yp^k_2(y,s)|^\frac 32ds\leq c_3\lambda_k^\frac 52.\ee The first component of the pressure satisfies
the equation
$$\De_y p^k_1(y,s)=-\div_y\div_y (\chi_{B(-x^k/\lambda_k, 1/\lambda_k)}(y)u^k(y,s)\otimes u^k(y,s)),\qquad y\in \mathbb R^3,$$
 for all possible values of $s$. For such a function, we have the standard estimate
\be\la{b18}\|p^k_1(\cdot,s)\|_{BMO(\mathbb R^3)}\leq c\ee
for all $s\in ]-(1-t_k)/\lambda^2_k,0[$.

We slightly change $p^k_1$  and  $p^k_2$ setting
$$\overline{p}_1^k(y,s)=p^k_1(y,s)-[p^k_1]_{B(1)}(s)\qquad\overline{p}_2^k(y,s)
=p^k_2(y,s)-[p^k_2]_{B(1)}(s)$$
so that $[\overline{p}_1^k]_{B(1)}(s)=0$ and $[\overline{p}_2^k]_{B(1)}(s)=0$.

Now, we pick up an arbitrary positive number $a$ and fix it. Then from (\ref{b17}) and (\ref{b18}) it follows that for sufficiently large $k$ we have
$$\int\limits_{Q(a)}|\overline{p}_1^k|^\frac 32de+\int\limits_{Q(a)}|\overline{p}_2^k|^\frac 32de\leq c_4(c_2,c_3,a).$$  Using the same bootstrap arguments, we can show that the following estimate is valid:
$$\|u^k\|_{C^\al(\overline{Q}(a/2)}\leq c_5(c_2,c_3,c_4,a)$$ for some positive number $\al<1/3$.  Indeed, the norm $\|u^k\|_{C^\al(\overline{Q}(a/2))}$ is estimated with the help of norms $\|u^k\|_{L_\infty( Q(a)))}$ and $\|\overline{p}^k\|_{L_\frac 32(Q(a))}$, where $\overline{p}^k=\overline{p}^k_1+\overline{p}^k_2$. Hence, using the diagonal Cantor procedure,
we can select subsequences such that for some positive $\al$ and for any positive $a$
$$u^k\rightarrow u\qquad \mbox{in}\,\, C^\al(Q(a)),$$
$$\overline{p}_1^k\rightharpoonup \overline{p}_1,\qquad \mbox{in}\,\,L_\frac 32(Q(a)),\qquad[\overline{p}_1]_{B(1)}(s)=0,$$
$$\overline{p}^k_2\rightharpoonup \overline{p}_2\qquad \mbox{in}\,\,L_\frac 32(Q(a)),\qquad[\overline{p}_2]_{B(1)}(s)=0.$$
Moreover,  $u$ is a bounded ancient solution with the total pressure $\overline{p}=\overline{p}_1+\overline{p}_2$, where $\overline{p}_1=p_{u\otimes u}$.

Next, for sufficiently large $k$, we get from (\ref{b17}) that
$$\sup\limits_{y\in \overline{B}(a)}\int\limits_{Q(a)}|\na p^k_2(y,s)|^\frac 32 ds\leq c_3\lambda^{\frac 52}_k.$$ Hence, $\na p_2=0$ in $Q(a)$ for any $a>0$.
So, $p_2(y,s)$ is identically zero. This allows us to conclude that the pair $u$ and $p_{u\otimes u}$ is a solution to the Navier-Stokes equations in the sense of distributions and thus $u$ is a nontrivial mild bounded ancient solution satisfying the condition
$|u(0,0)|=1$. Theorem \ref{rt8} is proved.

\setcounter{equation}{0}
\section{Axially Symmetric Suitable Weak \\ Solutions}Without loss of generality, the problem of local regularity of weak solutions
(not necessary being axially symmetric) to the Navier-Stokes equations  can be formulated as follows. Let us consider a pair of functions
$v\in L_3(Q)$
%$v\in W^{1,0}_2(Q)$
and
$q\in L_\frac 32(Q)$, defined in the unit space-time cylinder $Q=\mathcal C\times ]-1,0[$, where $\mathcal C= \{ x=(x',x_3),\, x'=(x_1,x_2),\, |x'|<1,\,|x_3|<1\}$ is the unit  spatial cylinder of $\mathbb R^3$, which satisfies the Navier-Stokes system
in $Q$ in the sense of distributions. The question we are interested in is under what additional conditions on $v$ and $q$, the space-time origin $z=(x,t)=0$ is a regular point of $v$. By the definition, the velocity $v$ is regular at the point $z=0$ if there exists a positive number $r$ such that $v$ is essentially bounded
in the space-time cylinder $Q(r)$. Here, $Q(r)=\mathcal C(r)\times ]-r^2,0[$ and
$\mathcal C(r)=\{|x'|<r,\,|x_3|<r\}$.
In contrast to traditional setting, we replace the usual balls with cylinders,
which is quite convenient in the case of axial symmetry. As usual, we set
$$\overline{v}=v_\varrho e_\varrho+v_3e_3\qquad\widehat{v}=v_\varphi e_\varphi$$
for $v=v_\varrho e_\varrho+v_\varphi e_\varphi+v_3e_3$.

We reformulate   our main results for these canonical domains. The general case is obtained by re-scaling.
\begin{theorem}\la{rt1} Assume that functions $v\in L_3(Q)$
%$v\in W^{1,0}_2(Q)$
and
$q\in L_\frac 32(Q)$ are an axially symmetric weak solution to the Navier-Stokes equations in $Q$. Let, in addition, %\be\la{r2}v\in L_\infty(\mathcal C\times %]-1,-a^2[)\ee
%for each $0<a<1$ and,
for some positive constant $C$,
\be\la{r3}|\overline{v}(x,t)|\leq \frac C{\sqrt{-t}}\ee
for almost all points $z=(x,t)\in Q$. Then $z=0$ is a regular point of $v$. \end{theorem}
\begin{theorem}\la{rt2} Assume that functions $v\in L_3(Q)$
%$v\in W^{1,0}_2(Q)$
and
$q\in L_\frac 32(Q)$ are an axially symmetric weak solution to the Navier-Stokes equations in $Q$. Let, in addition,  \be\la{r2}v\in L_\infty(\mathcal C\times ]-1,-a^2[)\ee
for each $0<a<1$ and
%(\ref{r2})
%condition (\ref{r2}) holds and,
\be\la{r4}|\overline{v}(x,t)|\leq \frac C{|x'|}\ee
for almost all points $z=(x,t)\in Q$ with some positive constant $C$. Then $z=0$ is a regular point of $v$.\end{theorem}

It is well-known due to Caffarelli-Kohn-Nirenberg that if
$z=(x,t)$ is singular (i.e., not regular) point of $v$, then there must be $x'=0$.
In other words, all singular points must belong to the axis of symmetry which is axis $x_3$.
\begin{lemma}\la{asl2} Assume that functions $v\in L_3(Q)$
%$v\in W^{1,0}_2(Q)$
and
$q\in L_\frac 32(Q)$ are an axially symmetric weak solution to the Navier-Stokes equations in $Q$. Let, in addition,
%condition  \be\la{r2}v\in L_\infty(\mathcal C\times ]-1,-a^2[)\ee
%for each $0<a<1$
condition (\ref{r2})
hold. Then following estimate is valid:
\be\la{as1}\ess \sup\limits_{z\in Q(1/2)}|\varrho v_\varphi(z)|\leq C(M)\Big(\int\limits_{Q(3/4)}|\varrho v_\varphi|^\frac {10}3dz\Big)^\frac 3{10},
\ee where
$$M=
%\frac 1{\sqrt{R}}
\Big(\int\limits_{Q(3/4)}|\overline{v}|^\frac {10}3dz\Big)^\frac 3{10}+1.$$
\end{lemma}
For the reader convenience, we put the proof of Lemma \ref{asl2} in Appendix II,
see also \ci{CL} and \ci{CSTY}. Here, we would like to notice the following.
\begin{remark}\la{asr3} Under the assumptions of Lemma \ref{asl2}, the pair $v$ and $q$ forms a suitable weak solution to the Navier-Stokes equations in $Q$. Hence, the right hand side of (\ref{as1}) is bounded from above. \end{remark}

We recall  that the Navier-Stokes equations are invariant with respect to the following scaling:
$$u(x,t)=\lb v(\lb x, \lambda^2 t), \qquad p(x,t)=\lb^2 q(\lb x, \lambda^2 t)$$
So, new functions $u$ and $p$ satisfy the Navier-Stokes equations in a suitable domain.

%It is also well understood that scaling invariant functionals play very important
%role in the regularity theory of the Navier-Stokes equations. To describe some of %them, we need additional notation:
With some additional notation
$$\mathcal C(x_0,R)=\{x\in \mathbb R^3\,\,\|\,\, x=(x',x_3),
\,\,x'=(x_1,x_2),\,\,$$$$|x'-x'_0|<R,\,\,|x_3-x_{03}|<R\},\qquad\mathcal
C(R)=\mathcal C(0,R),\qquad \mathcal C=\mathcal C(1);
$$$$ z=(x,t),\qquad z_0=(x_0,t_0),\qquad Q(z_0,R)=\mathcal C(x_0,R)
\times ]t_0-R^2,t_0[, $$$$
  Q(R)=Q(0,R),\qquad Q=Q(1),$$
we introduce certain scale-invariant functionals:
%Here is the list of some scaling invariant functionals which are going to be %exploited in what follows:
$$A(z_0,r;v)=\ess\sup\limits_{t_0-r^2<t<t_0}\frac 1r\int\limits_{\mathcal C(x_0,r)}|v(x,t)|^2dx,$$$$ E(z_0,r;v)=\frac 1r\int\limits_{Q(z_0,r)}|\na v|^2dz,\qquad D(z_0,r;q)=\frac 1{r^2}\int\limits_{Q(z_0,r)}|q|^\frac 32dz,$$
$$C(z_0,r;v)=\frac 1{r^2}\int\limits_{Q(z_0,r)}| v|^3dz,\qquad H(z_0,r;v)=\frac 1{r^3}\int\limits_{Q(z_0,r)}| v|^2dz,$$ $$M_{s,l}(z_0,r;v)=\frac 1{r^\kappa}\int\limits_{t_0-r^2}^{t_0}\Big(\int\limits_{\mathcal C(x_0,r)}|v|^sdx\Big)^\frac ls dt,$$
where $\kappa=l(\frac 3s+\frac 2l-1)$ and $s\geq1$, $l\geq 1$. As it was shown in \ci{S10}, the following inequality holds
$$C(z_0,r;f)\leq c A^\mu(z_0,r;f)(M_{s,l}(z_0,r;f))^\frac 1m(E(z_0,r;f)+$$\be\la{as2}+H(z_0,r;f))^\frac {m-1}m,\ee
 where
 $$\mu=\frac lm\Big(\frac 3s+\frac 3l-2\Big),\qquad m=2l \Big(\frac 3s+\frac 2l-\frac 32\Big),$$
 provided
 \be\la{as3}\frac 3s+\frac 2l-\frac 32\geq\max\Big\{ \frac 12-\frac 1s,\frac 1s-\frac 16\Big\}.\ee
 %It is easy to verify that condition (\ref{as3}) implies two facts: $m>1$ and $$0<\mu<1$.
 Actually, inequality (\ref{as2}) is but the result of application of H\"older's inequality and special Galiardo-Nireberg's inequality.

 The essential technical part of the proof of Theorem \ref{rt1} is the following lemma.
 \begin{lemma}\la{asl4}Under assumptions of Theorem \ref{rt1}, we have the estimate
 \be\la{as4}A(z_b,r;v)+E(z_b,r;v)+C(z_b,r;v)+D(z_b,r;q)\leq C_1<+\infty\ee
 for all $z_b$ and for all $r$ satisfying conditions
 \be\la{as5} z_b=(be_3,0), \quad b\in\mathbb R,\quad |b|\leq \frac 14,\qquad 0<r<\frac 14.\ee
 A constant $C_1$ depends only on the constant $C$ in (\ref{r3}),
 $\|v\|_{L_3(Q)}$, and \\ $\|q\|_{L_\frac 32(Q)}$.
 \end{lemma}
 \textsc{Proof} By Lemma \ref{asl2} and by Remark \ref{asr3}, we  have two inequalities:
 \be\la{as6}A(0,3/4;v)+E(0,3/4;v)\leq C_2<+\infty,\ee
\be\la{as7}|x'||v_\varphi(x,t)|\leq C_2\quad\mbox{for a.a.}\,  z=(x,t)\in Q(1/2).\ee
 Constant $C_2$ depends on the same arguments as constant $C_1$.

It follows from (\ref{as2}), that, for $s_1=\frac 74$ and $l_1=10$, inequality (\ref{as7}) takes the form
$$m_1= \frac {58}7,\quad \mu_1=\frac 1{58},\qquad M_{s_1,l_1}(z_b,r;\widehat{v})\leq cC_2^{10},$$
\be\la {as8}C(z_b,r;\widehat{v})\leq c A^\frac 1{58}(z_b,r;v)(C^{10}_2)^\frac 7{58}(E(z_b,r;v)+A(z_b,r;v))^\frac {51}{58}\ee provided
%for all $z_b$ and $r$ satisfying
conditions (\ref{as5}) hold.

To treat $\overline{v}$ which is the other part of the velocity $v$, we chose numbers $s_2=4$ and $l_2=\frac {12}7$. Then, for the same reasons as above, we find
$$m_2=\frac {10}7,\quad \mu_2 =\frac 3{14},\quad M_{s_2,l_2}(z_b,r;\overline{v})\leq c C^\frac {12}7,$$
\be\la {as9}C(z_b,r;\overline{v})\leq c A^\frac 3{14}(z_b,r;v)(C^\frac {12}7)^\frac 7{10}(E(z_b,r;v)+A(z_b,r;v))^\frac {3}{10}\ee for all $z_b$ and $r$ satisfying
conditions (\ref{as5}).

Adding (\ref{as8}) and (\ref{as9}), we show
$$C(z_b,r;v)\leq c\Big(C(z_b,r;\overline{v})+C(z_b,r;\widehat{v})\Big)\leq$$
\be\la{as10}\leq c\Big(A^\frac 1{58}(z_b,r;v)C_2^\frac {35}{29}(E(z_b,r;v)+A(z_b,r;v))^\frac {51}{58}+\ee
$$+A^\frac 3{14}(z_b,r;v)C^\frac {6}5(E(z_b,r;v)+A(z_b,r;v))^\frac {3}{10}\Big)$$
for the same  $z_b$ and $r$ as above.
%satisfying
%conditions (\ref{as5}).

Applying Young's inequality in (\ref{as10}), we arrive at the important estimate
\be\la{as11}C(z_b,r;v)\leq \varepsilon(E(z_b,r;v)+A(z_b,r;v))+f_1(\varepsilon, C, C_2),\ee
provided
%being valid for all $z_b$ and $r$ satisfying
conditions (\ref{as5}) hold. In (\ref{as11}), the positive number $\varepsilon$ is a parameter to pick up later. The rest of the proof is routine. In addition to (\ref{as11}), we consider
the local energy inequality
\be\la{as12}E(z_b,r/2;v)+A(z_b,r/2;v)\leq c\Big(C^\frac 23(z_b,r;v)+C(z_b,r;v)+D(z_b,r;q)\Big)\ee
%where  $z_b$ and $r$
%conditions (\ref{as5}),
and the decay estimate for the pressure field
\be\la{as13}D(z_b,\varrho;q)\leq c\Big[\frac \varrho r D(z_b,r;q)+\Big(\frac r\varrho\Big)^2C(z_b,r;v)\Big].\ee
Here,  $z_b$ and $r$ satisfy
conditions (\ref{as5}) and  $0<\varrho\leq  r$.  If we let
$$\mathcal E(r)=E(z_b,r;v)+A(z_b,r;v)+D(z_b,r;q),$$
then, for a fixed small positive number $\vartheta$, one can derive from (\ref{as12}) and (\ref{as13}) the following estimate
$$\mathcal E(\vartheta r)\leq c\Big(C^\frac 23(z_b,2\vartheta r;v)+C(z_b,2\vartheta r;v)+D(z_b,2\vartheta r;q)+$$$$+\vartheta D(z_b,r;q)+\frac 1{\vartheta^2}C(z_b,r;v)\Big)\leq$$
$$\leq c\Big[\vartheta D(z_b,r;q)+\frac 1{\vartheta^2}C(z_b,r;v)+\frac 1{\vartheta^\frac 43}C^\frac 23(z_b,r;v)\Big].$$
Now, the last two terms on the right hand side of the latter inequality can be majorized  with the help of  (\ref{as11}). As a result, we have
$$\mathcal E(\vartheta r)\leq c \Big(\vartheta+\frac\varepsilon{\vartheta^2}\Big)\mathcal E(r) +f_2(\varepsilon,\vartheta,C,C_2).$$
We first chose $\vartheta$ so that $c\vartheta<\frac 14$, pick up
$\varepsilon$ to provide the inequality $\frac {c\varepsilon}{\vartheta^2}<\frac 14$, and then we find
$$\mathcal E(\vartheta r)\leq \frac 12\mathcal E(r)+f_3(C,C_2).$$
The latter inequality can be easily iterated. After simple calculations, we derive
the relation
$$E(z_b,r;v)+A(z_b,r;v)+D(z_b,r;q)\leq c\Big(A(0,1/2;v)+E(0,1/2;v)+$$$$+D(0,1/2;q)+f_3(C,C_2)\Big)$$
with $z_b$ and $r$ satisfying
conditions (\ref{as5}). Lemma \ref{asl4} is proved.

To prove Theorem \ref{rt2}, we need an analogue of Lemma \ref{asl4}. Here, it is.
\begin{lemma}\la{asl5}Under assumptions of Theorem \ref{rt2}, estimate (\ref{as4}) is valid as well with constant $C_1$ depending only on the constant $C$ in (\ref{r4}),
 $\|v\|_{L_3(Q)}$, and  $\|q\|_{L_\frac 32(Q)}$.
%all the statements  of Lem-\\ma \ref{asl4} remain to be valid.
 \end{lemma}
 Lemma \ref{asl5} is proved in the same  way as  Lemma \ref{asl4} and even easier because  main inequality (\ref{as11}) can be established with the help  of the  case  $s=s_1$, $l=l_1$ only.

As it follows from conditions of Theorem \ref{rt2} and the statement of Lemma \ref{asl2}, the module of the velocity field grows not faster than $C/|x'|$ as $|x'|\to 0$. Moreover, the corresponding estimate is uniform in time. However, it turns out to be true under conditions of Theorem \ref{rt1} as well. More precisely, we have the following.
\begin{pro}\la{asp6} Assume that all conditions of Theorem \ref{rt1} hold.
Then
\be\la{as14}|v(x,t)|\leq \frac {C_1} {|x'|}\ee
for all $z=(x,t)\in Q(1/8)$.
%$x=(x',x_3)$ such that $|x'|<\frac 18$ $and |x_3|<\frac 18$
%and for all $t$ such that $-\Big(\frac 18\Big)^2<t<0$.
A constant $C_1$ depends only on the constant $C$ in (\ref{r3}),
 $\|v\|_{L_3(Q)}$, and  $\|q\|_{L_\frac 32(Q)}$.\end{pro}
\textsc{Proof} In view
%Having statement of Lemma
of (\ref{asl4}), we can  argue essentially as in
\ci{S11}.

Let us fix a point $x_0\in \mathcal C(1/8)$ and put
$r_0=|x'_0|$, $b_0=x_{03}$. So, we have $r_0<\frac 18$ and $|b_0|< \frac 18$.
Further, we introduce the following cylinders:
$$\mathcal P^1_{r_0}=\{r_0<|x'|<2r_0,\,|x_3|<r_0\},\quad\mathcal P^2_{r_0}=\{r_0/4<|x'|<3r_0,\,|x_3|<2r_0\}.$$
%With help of them, we can define other cylinders:
$$\mathcal P^1_{r_0}(b_0)=\mathcal P^1_{r_0}+b_0e_3,\quad\mathcal P^2_{r_0}(b_0)=\mathcal P^2_{r_0}+b_0e_3,$$
$$Q^1_{r_0}(b_0)=\mathcal P^1_{r_0}(b_0)\times ]-r^2_0,0[,\quad Q^2_{r_0}(b_0)=\mathcal P^2_{r_0}(b_0)\times ]-(2r_0)^2,0[.$$
Now, let us scale our functions so that
$$x=r_0y+b_0e_3, \quad t=r_0^2s, \quad u(y,s)=r_0v(x,t),\quad p(y,s)=r_0^2q(x,t).$$
As it was shown in \ci{S11}, there exists a continuous nondecreasing function $\Phi_:\mathbb R_+\to\mathbb R_+$, $\mathbb R_+=\{s>0\}$, such that
$$\sup\limits_{(y,s)\in Q^1_1(0)}|u(u,s)|+|\na u(y,s)|\leq\Phi\Big(\sup\limits_{-2^2<s<0}\int\limits_{\mathcal P^2_1(0)}|u(y,s)|^2dy$$
\be\la{as15} +\int\limits_{Q^2_1(0)}|\na u|^2dy\,ds+
\int\limits_{Q^2_1(0)}| u|^3dy\,ds+\int\limits_{Q^2_1(0)}|p|^\frac 32dy\,ds\Big).\ee
After making inverse scaling in (\ref{as15}), we find
$$\sup\limits_{z\in Q^1_{r_0}(b_0)}r_0|u(x,t)|+r^2_0|\na u(x,t)|\leq \Phi\Big(cA(z_{b_0},3r_0;v)+cE(z_{b_0},3r_0;v)+$$$$+cC(z_{b_0},3r_0;v)+
cD(z_{b_0},3r_0;q)\Big)\leq \Phi\Big(4cC_1\Big).$$ It remains to apply Lemma \ref{asl4} and complete
the proof of the proposition. Proposition \ref{asp6} is proved.

\setcounter{equation}{0}
\section{Proof of Theorems \ref{rt1} and \ref{rt2} }
Using  Lemmata \ref{asl2}, \ref{asl4}, \ref{asl5}, Remark \ref{asr3},
Proposition \ref{asp6} and scaling arguments, we may assume (without loss of generality) that our solution
$v$ and $q$  have the following properties:
\be\la{p1}\sup\limits_{0<r\leq 1} \Big(A(0,r;v)+E(0,r;v)+C(0,r;v)+D(0,r;q)\Big)=A_1<+\infty,\ee
\be\la{p2}\ess \sup\limits_{z=(x,t)\in Q}|x'||v(x,t)|=A_2<+\infty.\ee
We  may also assume that the function $v$ is H\"older continuous in the completion
of the set $\mathcal C\times ]-1,-a^2[$ for any $0<a<1$.

Introducing  functions
$$H(t)=\sup\limits_{x\in \mathcal C}|v(x,t)|,\quad h(t)=\sup\limits_{-1<\tau\leq t}H(t),$$
 let us suppose that our statement is wrong, i.e., $z=0$ is a singular point. Then there are  sequences $x_k\in \overline{\mathcal C}$ and $-1<t_k<0$,
having the following properties:
$$h(t_k)=H(t_k)=M_k=|v(x_k,t_k)|\to +\infty \qquad \mbox{as}\quad k\to+\infty.$$
We scale our functions $v$ and $q$ so that scaled functions possess axial symmetry:
$$u^k(y,s)=\lambda_kv(\lambda_k y',x_{3k}+\lambda_k y_3, t_k+\lambda^2_ks),\qquad
\lambda_k =\frac  1{M_k},$$$$p^k(y,s)=\lambda^2_kq(\lambda_k y',x_{3k}+\lambda_k y_3, t_k+\lambda^2_ks).$$
These functions satisfy the Navier-Stokes equations in $Q(M_k)$. Moreover,
\be\la{p3} |u^k(y'_k,0,0)|=1,\qquad y'_k=M_kx_k'.\ee
According to (\ref{p2}),
%we have
$$|y'_k|\leq A_2$$
for all $k\in \mathbb N$. Thus, without loss of generality, we may assume that
\be\la{p4}y'_k\to y'_*\qquad \mbox{as}\quad k\to +\infty.\ee

Now, let us see what happens as $k\to +\infty$. By the identity
\be\la{p5}\sup\limits_{e=(y,s)\in \mathcal C(M_k)}|u^k(e)|=1\ee and by (\ref{p1}), we can select subsequences (still
denote as the entire sequence) such that
\be\la{p6} u^k{\stackrel{\star}{\rightharpoonup}}\,u\qquad\mbox {in}\qquad
L_{\infty}(Q(a)),\ee
and
\be\la{p7}p^k\rightharpoonup\, p \qquad\mbox {in}\qquad
L_{\frac 32}(Q(a))\ee
for any $a>0$. Functions $u$ and $p$ are defined on $Q_-=\mathbb R^3\times ]-\infty,0[$. Obviously, they possess the following properties:
\be\la{p8}\ess\sup\limits_{e\in Q_-}|u(e)|\leq 1,\ee
\be\la{p9}\sup\limits_{0<r<+\infty } \Big(A(0,r;u)+E(0,r;u)+C(0,r;u)+D(0,r;p)\Big)\leq A_1,\ee
\be\la{p10}\ess \sup\limits_{e=(y,s)\in Q_-}|y'||u(y,s)|\leq A_2.\ee

Now, our aim is to show that $u$ and $p$ satisfy the Navier-Stokes equations $Q_-$
and $u$ is smooth enough to obey the identity
\be\la{p11}|u(y'_*,0,0)|=1.\ee
To this end, we fix an arbitrary positive number $a>0$ and consider numbers $k$ so big that $a<M_k/4$. We know that $u^k$ satisfies the nonhomogeneous heat equation
of the form
$$\pa_t u^k-\De u^k=-\div F^k\qquad \mbox{in}\quad Q(4a),
$$
where $F^k=u^k\otimes u^k +p^k \mathbb I$ and
%we know that
$$\|F^k\|_{\frac 32,Q(4a)}\leq c_1(a)<\infty.$$
This is implies the following fact, see \ci{LSU},
$$\|\na u^k\|_{\frac 32,Q(3a)}\leq c_2(a)<\infty.$$
Now, we can interpret the pair $u^k$ and $p^k$ as a solution to the nonhomogeneous Stokes system
\be\la{p12}\pa_t u^k-\De u^k+\na p^k= f^k,\quad\div u^k=0\qquad \mbox{in}\quad Q(3a),\ee
where $f^k=-u^k\cdot \na u^k$ is the right hand side having the property
%and we have the estimate
$$\|f^k\|_{\frac 32,Q(3a)}\leq c_2(a). $$
Then, according to the local regularity theory for the Stokes system, see \ci{S8}, we can state that
$$\|\pa_tu^k\|_{\frac 32,Q(2a)}+\|\na^2u^k\|_{\frac 32,Q(2a)}+\|\na k^k\|_{\frac 32,Q(2a)}\leq c_3(a).$$
The latter, together with the embedding theorem, implies
%we have
$$\|\na u^k\|_{3,\frac 32,(Q(2a))}+\|p^k\|_{3,\frac 32,Q(2a)}\leq c_4(a).
$$
In turn, this improves integrability of the right hand side in (\ref{p12})
$$\|f^k\|_{3,\frac 32,Q(2a)}\leq c_4(a).$$
Therefore, by the  local regularity theory,
$$\|\pa_tu^k\|_{3,\frac 32,Q(2a)}+\|\na^2u^k\|_{3,\frac 32,Q(2a)}+\|\na k^k\|_{3,\frac 32,Q(2)}\leq c_5(a).$$ Applying the imbedding theorem once more,
we find
$$\|\na u^k\|_{6,\frac 32,Q(2a)}+\|p^k\|_{6,\frac 32,Q(2a)}\leq c_6(a).
$$
The local regularity theory leads then to the estimate
$$\|\pa_tu^k\|_{6,\frac 32,Q(a)}+\|\na^2u^k\|_{6,\frac 32,Q(a)}+\|\na p^k\|_{6,\frac 32,Q(a)}\leq c_7(a).$$
By the embedding theorem, sequence $u^k$ is uniformly bounded in the parabo-\\lic
H\"older space $C^\frac 12(\overline{Q}(a/2))$. Hence, without loss of generality, one may assume that
$$ u^k\rightarrow  u\qquad \mbox{in}\quad C^\frac 14(\overline{Q}(a/2)). $$
This means that the pair $u$ and $p$ obeys the Navier-Stokes system and (\ref{p11})
holds. So, the function $u$ is the so-called bounded ancient solution to the Navier-Stokes system
which is, in addition, axially symmetric and satisfies the decay estimate (\ref{p10}).
As it was shown in \ci{KNSS}, such a solution must be identically zero. But this contradicts (\ref{p11}). Theorems (\ref{rt1}) and (\ref{rt2}) are proved.

\setcounter{equation}{0}
\section{Appendix I: Proof of Theorem \ref{bt2} }

In what follows, we need a few known regularity results.
%Now, let us formulate several auxiliary lemmata.
%\begin{lemma}\la{pl1} For any $F=L_\infty(\mathbb R^n;\mathbb M^{n\times n})$, there %exists a unique function $q_F\in BMO(\mathbb R^n)$ such that $[q_F]_{B(1)}=0$ and %$$\De q_F=-\div \div F=F_{ij,ij} \qquad \mbox{in}\quad \mathbb R^3$$
%in the sense of distributions. Moreover, function $q_F$ meets the estimate
%$$\|q_F\|_{BMO(\mathbb R^n)}\leq c(n)\|F\|_{\infty,\mathbb R^n}.$$\end{lemma}
%To state Lemma \ref{pl1}, the following notation has been used.
%$[f]_\Om$
%is the mean value of a function $f$ over a spatial domain $\Om\in\mathbb R^n$. In %turn, the mean value of a function $g$ over a space-time domain $Q$ is denoted by %$(g)_Q$.
%Comma in subscripts stands for differentiation with respect to the corresponding %spatial variable, for example, $f_{,i}=\pa f/\pa x_i$. Summation over repeated Latin %indices
%running from 1 to $n$ is adopted.
\begin{lemma}\la{pl2} Assume that functions $f\in L_m(B(2))$ and $q\in L_m(B(2))$ satisfy the equation
$$\De q=-\div f\qquad \mbox{in}\quad B(2).$$
Then
$$\int\limits_{B(1)}|\na q|^mdx\leq c(m,n)\Big(\int\limits_{B(2)}|f|^mdx+\int\limits_{B(2)}|q-[q]_{B(2)}|^mdx\Big).
$$\end{lemma}
\begin{lemma}\la{pl3} Assume that functions $f\in L_m(Q(2))$ and $u\in W^{1,0}_m(Q(2))$ satisfy the equation
$$\pa_t u-\De u=f\qquad \mbox{in}\quad Q(2).$$
Then $u\in W^{2,1}_m(Q(1))$ and the following estimate is valid:
$$\|\pa_t u\|_{m,Q(1)}+\|\na^2 u\|_{m,Q(1)}\leq c(m,n)\Big[\|f\|_{m,Q(2)}+\|u\|_{W^{1,0}_m(Q(2))}\Big].$$\end{lemma}
 %Lemma \ref{pl1} is proved with the help of the singular integral theory,
% see \ci{St}.
Proof of Lemmata \ref{pl2} and \ref{pl3} can be found, for example, in \ci{LU} and \ci{LSU}.

\textsc{Proof of Theorem \ref{bt2}}: \textsc{Step 1}.\textsc{Energy estimate}. Take
an arbitrary number $t_0<0$ and fix it. Let $k_\varepsilon(z)$ be a standard smoothing kernel and let
%. We use the following notation for mollified functions:
$$F^\varepsilon(z)=\int\limits_{Q_-}k_\varepsilon(z-z')F(z')dz',\qquad F=u\otimes u,$$
$$u^\varepsilon(z)=\int\limits_{Q_-}k_\varepsilon(z-z')u(z')dz'.$$

Assume that $w\in {\stackrel{\circ}{C}}{^\infty_{0}}(Q^{t_0}_-)$, where $Q^{t_0}_-=\mathbb R^n\times ]-\infty,t_0[$. Obviously, $w^\varepsilon\in {\stackrel{\circ}{C}}{^\infty_{0}}(Q_-)$ for sufficiently small $\varepsilon$. %$(0<\varepsilon<\varepsilon(t_0))$, we have . Then
Using known properties of smoothing kernel and Definition \ref{bd1}, we find
$$\int\limits_{Q_-}w\cdot (\pa_t u^\varepsilon+\div F^\varepsilon-\De u^\varepsilon)dz=0,\qquad \forall w\in {\stackrel{\circ}{C}}{^\infty_{0}}(Q^{t_0}_-).$$
There exists a smooth function $p_\varepsilon$
with the following property
\be\la{b1}\pa_tu^\varepsilon+\div F^\varepsilon-\De u^\varepsilon=-\na p_\varepsilon,\qquad \div u^\varepsilon=0\ee
in $Q_-^{t_0}$. Splitting  pressure $p_\varepsilon$ into two parts
\be\la{b2} p_\varepsilon=p_{F^\varepsilon}+\widetilde{p_\varepsilon}.\ee
and observing
%It is not difficult to show
that the function $\na p_{F^\varepsilon}$ is bounded in $Q_-^{t_0}$, one can conclude that, by
% So, it follows from
(\ref{b1}) and (\ref{b2}),
% that
$$\De \widetilde{p_\varepsilon}=0 \qquad \mbox{in}\quad Q^{t_0}_-,\qquad
\na \widetilde{p_\varepsilon}\in L_\infty(Q^{t_0}_-;\mathbb R^n).$$
According to Liouville's theorem for harmonic   functions,
$$ \na \widetilde{p_\varepsilon}(x,t)=a_\varepsilon(t), \qquad x\in\mathbb R^n,\quad
 -\infty<t\leq t_0.$$ So, we have
\be\la{b3}\pa_tu^\varepsilon+\div F^\varepsilon-\De u^\varepsilon=-\na p_{F^\varepsilon}-a_\varepsilon,\qquad \div u^\varepsilon=0\ee
in $Q_-^{t_0}$.

Now, let us introduce new auxiliary functions
$$b_{\varepsilon t_0}(t)=\int\limits^{t}_{t_0-1}a_\varepsilon(\tau)d\tau,\qquad t_0-1\leq t\leq t_0,$$
$$v_\varepsilon(x,t)=u^\varepsilon(x,t)+b_{\varepsilon t_0}(t),\qquad z=(x,t)\in Q^{t_0}.$$
Using them, one may reduce system (\ref{b3}) to the form
\be\la{b4}\pa_tv_\varepsilon-\De v_\varepsilon=-\div F^\varepsilon-\na p_{F^\varepsilon},\qquad \div v_\varepsilon=0\ee
in $Q_-^{t_0}$.

Let $\varphi_{x_0}(x)=\varphi(x-x_0)$ for  a fixed cut-off function $\varphi$ satisfying  the conditions
$$0\leq\varphi\leq 1, \qquad \varphi\equiv 1\quad \mbox{in}\quad B(1),\qquad
\supp \varphi\subset B(2).$$

%Now, we can
To derive the energy identity, let us  multiply (\ref{b4}) by
$\varphi^2_{x_0}v_\varepsilon$ and integrate the product by parts. As a result, we have
$$I(t)=\int\limits_{\mathbb R^n}\varphi^2_{x_0}(x)|v_\varepsilon(x,t)|^2dx
+2\int\limits^{t}_{t_0-1}\int\limits_{\mathbb R^n}\varphi^2_{x_0}|\na v_\varepsilon|^2dxdt'=$$
$$=\int\limits_{\mathbb R^n}\varphi^2_{x_0}(x)|v_\varepsilon(x,t_0-1)|^2dx+
\int\limits^{t}_{t_0-1}\int\limits_{\mathbb R^n}\De \varphi^2_{x_0}| v_\varepsilon|^2dxdt'+$$
$$+\int\limits^{t}_{t_0-1}\int\limits_{\mathbb R^n}(p_{F^\varepsilon}-[p_{F^\varepsilon}]_{B(x_0,2)})v_\varepsilon\cdot \na \varphi^2_{x_0}dxdt'+$$$$+\int\limits^{t}_{t_0-1}\int\limits_{\mathbb R^n}(F^\varepsilon-[F^\varepsilon]_{B(x_0,2)}):\cdot \na (\varphi^2_{x_0}v_\varepsilon)dxdt'.$$

Introducing
$$\al_\varepsilon(t)=\sup\limits_{x_0\in\mathbb R^n}\int\limits_{B(x_0,1)}
|v_\varepsilon(x,t)|^2dx$$ and taking into account that $v_\varepsilon(\cdot,t_0-1)=u^\varepsilon(\cdot,t_0-1)$ and $|u^\varepsilon(\cdot,t_0-1)|\leq 1$, we can estimate the right hand side of the energy identity in the following way
$$I(t)\leq c(n)+c(n)\int\limits^{t}_{t_0-1}\al_\varepsilon(t')dt'
+$$$$+c(n)\Big(\int\limits^{t_0}_{t_0-1}\int\limits_{B(x_0,2)}|p_{F^\varepsilon}-[p_{F^\varepsilon}]_{B(x_0,2)}|^2dxdt\Big)^\frac 12
\Big(\int\limits^{t}_{t_0-1}\al_\varepsilon(t')dt'\Big)^\frac 12+$$
\be\la{b5}+c(n)\Big(\int\limits^{t_0}_{t_0-1}\int\limits_{B(x_0,2)}|F^\varepsilon-[F^\varepsilon]_{B(x_0,2)}|^2dxdt\Big)^\frac 12
\Big(\int\limits^{t}_{t_0-1}\int\limits_{\mathbb R^n}\varphi^2_{x_0}
|\na v_\varepsilon|^2dxdt'+\ee$$+\int\limits^{t}_{t_0-1}\al_\varepsilon(t')dt'\Big)^\frac 12, \qquad t_0-1\leq t\leq t_0.$$
Next, since $|F^\varepsilon|\leq c(n)$, we find two estimates
$$\int\limits^{t_0}_{t_0-1}\int\limits_{B(x_0,2)}|F^\varepsilon-[F^\varepsilon]_{B(x_0,2)}|^2dxdt\leq c(n)$$
and
$$\int\limits^{t_0}_{t_0-1}\int\limits_{B(x_0,2)}|p_{F^\varepsilon}-[p_{F^\varepsilon}]_{B(x_0,2)}|^2dxdt\leq c(n)\|p_{F^\varepsilon}\|^2_{L_\infty(-\infty,t_0;BMO(\mathbb R^n))}$$$$\leq c(n)\|F^\varepsilon\|^2_{L_\infty(Q^{t_0}_-)}\leq c(n).$$
The latter estimates, together with (\ref{b5}), imply the inequalities
$$\al_\varepsilon(t)\leq c(n)\Big(1+\int\limits^{t}_{t_0-1}\al_\varepsilon(t')dt'\Big), \qquad t_0-1\leq t\leq t_0$$
and
$$\sup\limits_{x_0\in \mathbb R^n}\int\limits^{t_0}_{t_0-1}\int\limits_{B(x_0,1)}
|\na v_\varepsilon|^2dxdt\leq c(n)\Big(1+\int\limits^{t_0}_{t_0-1}\al_\varepsilon(t)dt\Big).$$
Applying known arguments, we can conclude
\be\la{b6}\sup\limits_{t_0-1\leq t\leq t_0}\al_\varepsilon(t)+\sup\limits_{x_0\in \mathbb R^n}\int\limits^{t_0}_{t_0-1}\int\limits_{B(x_0,1)}
|\na u^\varepsilon|^2dxdt\leq c(n).\ee
It should be emphasized that the right  hand size in (\ref{b6}) is independent of $t_0$. In particular, estimate (\ref{b6}) allows to show
$$\sup\limits_{t_0-1\leq t\leq t_0}b_{\varepsilon t_0}(t)\leq c(n).$$

Now, let us see what happens if $\varepsilon\to 0$. Selecting a subsequence if necessary and taking the limit as $\varepsilon\to 0$, we state that:
%have the following facts:
$$b_{\varepsilon t_0}{\stackrel{\star}{\rightharpoonup}}\,b_{t_0}\qquad\mbox {in}\qquad L_\infty(t_0-1,t_0;\mathbb R^n);$$
the estimate
\be\la{b7}\|b_{t_0}\|_{L_\infty(t_0-1,t_0)}+\sup\limits_{x_0\in \mathbb R^n}\int\limits_{t_0-1}^{t_0}\int\limits_{B(x_0,1)}
|\na u|^2dxdt\leq c(n)<+\infty\ee
is valid for all $t_0<0$;  the system
$$\pa_tu^{t_0}+\div u\otimes u-\De u=-\na p_{u\otimes u},\qquad \div u=0$$
holds in $Q^{t_0}$ in the sense of distributions.

The case $t_0=0$ can be treated by passing to the limit as $t_0\to 0$.

\textsc{Step 2, Bootstrap Arguments} By (\ref{b7}),
$$f=\div F=u\cdot \na u \in {\cal L}_2(Q_-;\mathbb R^n).$$
Then Lemma \ref{pl2}, together with shifts, shows that
$$\na p_{u\otimes u}\in {\cal L}_2(Q_-;\mathbb R^n).$$
Next, obviously, the function $u^{t_0}$ satisfies the system of equations
$$\pa_tu^{t_0}-\De u^{t_0}=-u\cdot \na u-\na p_{u\otimes u}\in {\cal L}_2(Q_-;\mathbb R^n),$$
which allows us to apply  Lemma \ref{pl3}  and conclude that
$$u^{t_0}\in W^{2,1}_2(Q(z_0,\tau_2);\mathbb R^n), \qquad 1/2 <\tau_2<\tau_1=1.$$
Moreover, the estimate
$$\|u^{t_0}\|_{W^{2,1}_2(Q(z_0,\tau_2))}\leq c(n,\tau_2)$$
holds for any $z_0=(x_0,t_0)$, where $x_0\in \mathbb R^n$ and $t_0\leq 0$.
Applying  the parabolic embedding theorem, see \ci{LSU}, we can state that
$$\na u^{t_0}=\na u\in W^{1,0}_{m_2}(Q(z_0,\tau_2);\mathbb R^n),$$
where
$$\frac 1{m_2}=\frac 1{m_1}-\frac 1{n+2},\qquad m_1=2.$$
By Lemma \ref{pl2}, by shifts, and by scaling,
 $$\int\limits_{B(x_0,\tau_3')}|\na p_{u\otimes u}(\cdot,t)|^{m_2}dx\leq c(n,\tau_2,\tau_3')\Big[\int\limits_{B(x_0,\tau_3')}|\na  u(\cdot,t)|^{m_2}dx
 +1\Big]$$ for $1/2<\tau'_3<\tau_2$.
 In turn, Lemma \ref{pl3} provides two statements:
 $$u^{t_0}\in W^{2,1}_{m_2}(Q(z_0,\tau_3);\mathbb R^n), \qquad 1/2 <\tau_3<\tau_3'$$ and
 $$\|u^{t_0}\|_{W^{2,1}_{m_2}(Q(z_0,\tau_3))}\leq c(n,\tau_3,\tau_3').$$
 Then, again, by the embedding theorem, we find
 $$\na u^{t_0}=\na u\in W^{1,0}_{m_3}(Q(z_0,\tau_3);\mathbb R^n)$$
with
$$\frac 1{m_3}=\frac 1{m_2}-\frac 1{n+2}.$$

Now, let us take an arbitrary large number $m>2$ and fix it.  Find $\al$
as an unique solution to the equation
$$\frac 1m=\frac 12 -\frac \al{n+2}.$$ Next, for $k_0=[\al]+1$, where $[\al]$ is the entire part of the number $\al$,  determine the number $m_{k_0+1}$
satisfying the identity
$$\frac 1{ m_{k_0+1}}=\frac 12 -\frac {k_0}{n+2}.$$ Obviously, $m_{k_0+1}>m$.
Setting
$$\tau_{k+1}=\tau_k-\frac 14\frac 1{2^k}, \qquad \tau_1=1,\qquad k=1,2,,,$$
and repeating our previous arguments $k_0$ times, we conclude that:
$$u^{t_0}\in W^{2,1}_{{m_{k_0+1}}}(Q(z_0,\tau_{k_0+1});\mathbb R^n)$$ and
 $$\|u^{t_0}\|_{W^{2,1}_{m_{k_0+1}}(Q(z_0,\tau_{k_0+1}))}\leq c(n,m).$$
 Since $\tau_k>1/2$ for any natural numbers $k$, we complete the proof of Theorem
 \ref{bt2}. Theorem \ref{bt2} is proved.

 We can exclude the pressure field completely by considering the equations
for vorticity $\om=\na \wedge u$. In dimensions three, differentiability properties of $\om$ are described by the following theorem.
\begin{lemma}\la{bt5} Let $u$ be an arbitrary bounded ancient solution. For any $m>1$, we have
$$\om=\na\wedge u\in {\cal W}^{2,1}_m(Q_-;\mathbb R^3)$$
and
$$\pa_t\om +u\cdot \na \om -\De \om=\om\cdot\na u\qquad\mbox{a.e. in}\quad Q_-.$$ \end{lemma}
\begin{remark}\la{br6} We could continue investigations of regularity for solutions to the vorticity equations further and it would be a good exercise. However, regularity results stated in Theorem \ref{bt5} are sufficient for our purposes.\end{remark}
\begin{remark}\la{br7} %By the embedding theorems,
Functions $\om$ and $\na\om$ are H\"older continuous in $Q_-$ and their norms in H\"older spaces  are uniformly bounded there, see \ci{LSU}.\end{remark}

 \textsc{Proof of Lemma \ref{bt5}} Let us consider the case $n=3$. The case $n=2$ is in fact easier. So, we have
 $$\pa_t \om-\De \om=\om\cdot\na u-u\cdot \na\om\equiv f.$$
 Take an arbitrary number $m>2$ and fix it. By Theorem \ref{bt2},
 %the right hand side has the following property
 $$|f|\leq c(n)(|\na^2u|+|\na u|^2)\in L_m(Q(z_0,2))$$
 and the norm of $f$ in $L_m(Q(z_0,2))$ is dominated  by a constant depending only on
 $m$.
  %and being independent of $z_0$.
  It remains to apply Lemma \ref{pl3} and complete the proof of Lemma \ref{bt5}. Lemma \ref{bt5} is proved.

\setcounter{equation}{0}
\section{Appendix II: Proof of Lemma \ref{asl2} }

According to the local regularity theory of the Navier-Stokes equations,
see, for instance, \ci{S8}, one can easily show that the pair $v$ and $q$
has the following differentiability properties:
\be\la{a21}v\in W^{2,1}_\frac 32(\mathcal C(a)\times ]-a^2,-b^2[),\qquad
q\in  W^{1,0}_\frac 32(\mathcal C(a)\times ]-a^2,-b^2[)\ee
and
\be\la{a22} v\quad\mbox{is H\"older continuous in the completion of}\quad \mathcal C(a)\times ]-a^2,-b^2[\ee
for any $0<b\leq a<1$.

Now, we fix a number $m\geq 2$, multiply the equation for the velocity component $v_\varphi$
by $ru|u|^{m-2}$, where $u=rv_\varphi$,  and integrate the product by parts.
In  view of (\ref{a21}) and (\ref{a22}), we find the following identity
for $\om=|u|^m$
$$ \int\limits_{\mathcal C}\psi^2(x,t_*)|\om(x,t_*)|^2dx+\frac {2(m-1)}m\int\limits^{t_*}_{-1}\int\limits_{\mathcal C}\psi^2|\na \om(x,t)|^2dxdt$$
\be\la{a23}=\int\limits^{t_*}_{-1}\int\limits_{\mathcal C}| \om(x,t)|^2\Big(
\pa_t\psi^2+\overline{v}\cdot\na\psi^2+\De_2\psi^2+\frac{3\psi^2_{,r}}r\Big)dxdt.\ee
It is valid for all $-1<t_*<0$ and for all  cut-off functions vanishing
in a neighborhood of the boundary of the space-time cylinder $\mathcal C\times ]-1,1[$. Here,  $\De_2\psi^2=\psi^2_{,rr}+\psi^2_{,33}$. So, (\ref{a23})  means that the energy norm of $\psi\om$ is finite, i.e.,
$$ |\psi\om |^2_{2,\mathcal C\times ]-1,t_*[}
\equiv \ess\sup\limits_{ t\in ]-1,t_*[ }\int\limits_{\mathcal C} |\psi\om(x,t) |^2dx+\int\limits^{t_*}_{-1}\int\limits_{\mathcal C}|\na (\psi\om)|^2dxdt$$
\be\la{a24}\leq c\int\limits^{t_*}_{-1}\int\limits_{\mathcal C}| \om(x,t)|^2\Big(
\pa_t\psi^2+\overline{v}\cdot\na\psi^2+\De_2\psi^2+\frac{3\psi^2_{,r}}r
+|\na \psi|^2\Big)dxdt\ee
for any $-1<t_*<0$.

No, let us specify our cut-off function $\psi$ setting
$\psi(x,t)=\Phi(x)\chi(t)$ and assuming that new smooth functions $0\leq \Phi\leq 1$ and $0\leq \chi\leq 1$
meet the following properties:
$$\supp \Phi\in \mathcal C(r_1),\qquad \Phi\equiv 1\quad\mbox{in} \quad\mathcal C(r),$$
$$|\na \Phi|\leq \frac c{r_1-r},\qquad |\na^2 \Phi|\leq \frac c{(r_1-r)^2},\qquad |\pa_t\chi |\leq \frac c{(r_1-r)^2}.$$ Here, arbitrary fixed number $r$ and $r_1$
satisfy the condition
\be\la{a25}\frac 12<r<r_1<\frac 34.\ee
If we let $\widetilde{Q}=\mathcal C(r_1)\times ]-r_1^2,t_*[$, then
\be\la{a26}|\psi\om|^2_{2,\widetilde{Q}}\leq \frac c{(r_1-r)^2}\int\limits_{\widetilde{Q}}|\om|^2dz+c I,\ee
where
$$I=\frac 1{r_1-r}\int\limits_{\widetilde{Q}}|\psi\om||\om||\overline{v}|dz.$$
By H\"older's inequality,
$$I\leq \frac 1{r_1-r}\Big(\int\limits_{\widetilde{Q}}|\overline{v}|^\frac {10}3dz\Big)^\frac 3{10}\Big(\int\limits_{\widetilde{Q}}|\om|^\frac 52dz\Big)^\frac 25\Big(\int\limits_{\widetilde{Q}}|\psi\om|^\frac {10}3\Big)^\frac 3{10}.$$
The  left hand side of (\ref{a25}) can be evaluated from below
with the help of the well-known multiplicative inequality
\be\la{a27}\Big(\int\limits_{\widetilde{Q}}|\psi\om|^\frac {10}3dz\Big)^\frac 3{10}\leq
c|\psi\om|_{2,\widetilde{Q}}.\ee
Now, taking into account restriction (\ref{a25}), it is not so difficult to derive from (\ref{a26}) and (\ref{a27}) the following estimate
\be\la{a28}\Big(\int\limits_{\widetilde{Q}}|\psi\om|^\frac {10}3dz\Big)^\frac 3{10}\leq c M \frac {\sqrt{r_1}}{r_1-r}\Big(\int\limits_{\widetilde{Q}}|\om|^\frac 52dz\Big)^\frac 25.\ee
Setting
$$m=m_k=\Big(\frac 43\Big)^k,\quad r_1=r^{(k)}=\frac 12+\frac 1{2^{k+1}},\quad
r=r^{(k+1)},$$$$ \psi=\psi_k,\quad\widetilde{Q}_k=\mathcal C(r^{(k)})\times ]-(r^{(k)})^2,t_*[,\qquad k\in \mathbb N.$$
one  can  reduce (\ref{a28}) to the form
\be\la{a29}\Big(\int\limits_{\widetilde{Q}_k}|\psi_k u|^\frac {10m_k}3dz\Big)^\frac 3{10 m_k}\leq c M \frac {\sqrt{r^{(k)}}}{r^{(k)}-r^{(k+1)}}\Big(\int\limits_{\widetilde{Q}_k}|u|^\frac {5m_k}2dz\Big)^\frac 2{5m_k}\ee
The only difference with respect to the usual Moser's technique is that one should
take the limit as $t_*\to 0$ step-by-step. For example, for $k=1$, the integral
$$\int\limits_{Q(3/4)}| u|^\frac {10}3dz$$ is finite and therefore we can pass to the limit as $t_*\to 0$ in (\ref{a28}).  Then we may pass to the limit as $t_*\to 0$ in (\ref{a28}) for $k=2$ and so on. Tending $k\to +\infty$, we complete the proof of Lemma \ref{asl2} in more or less standard way. Lemma \ref{asl2} is proved.

%G. Seregin\\
%Steklov Institute of Mathematics at St.Petersburg, \\
%St.Petersburg, Russia
%\\
%\\

\end{document}